\documentclass{article}
\author{Oleg Pikhurko\thanks{Supported by a Research Fellowship, St.\ John's College, Cambridge.}\\
 DPMMS, Centre for Mathematical Sciences\\
 Cambridge University, Cambridge~CB3~0WB, England\\
 E-mail: {\tt O.Pikhurko@dpmms.cam.ac.uk}}

\usepackage{amssymb}
\newcommand{\eqref}[1]{\mbox{\rm(\ref{#1})}}
\newcommand{\pl}{}

\newcommand{\B}[1]{{\bf #1}}
\newcommand{\I}[1]{{\mathbb #1}}

\newcommand{\comment}[1]{}
\newcommand{\qed}{\nolinebreak\mbox{\hspace{5 true pt}%
\rule[-0.85 true pt]{3.9 true pt}{8.1 true pt}}}

\newtheorem{theorem}{Theorem}
\newtheorem{lemma}[theorem]{Lemma}
\newtheorem{corollary}[theorem]{Corollary}

\newcommand{\rcr}[1]{\ref{cr:\pl:#1}}

\begin{document}

\renewcommand{\B}[1]{\overline{#1}}

\title{Testing the Existence of a Supporting Plane}
\maketitle

\begin{abstract}
 We present an algorithm testing wheather, for given four vectors in
$\I R^3$, there is a plane through the origin such that all four
vectors fall into the same open halfspace.
 \end{abstract}

The following problem was communicated to me by Igor Komarov. Suppose
we are given four vectors $\B x_1,\dots,\B x_4\in\I R^3$. How to test
if there is a {\em supporting plane}, that is, a plane $P$ such that
all four vectors lie in the same open halfspace with respect to this
plane? (When we talk about planes or lines, we mean {\em linear
subspaces}, that is, containing the origin.)

This was needed for Komarov's survey paper~\cite{komarov:01}: the
carbon atoms in a molecule are classified there into two types
depending on the existence of a supporting plane (with respect to
the four vectors corresponding to the four ties at the carbon atom).

We did not find an explicit algorithm in the literature and, in order
to keep~\cite{komarov:01} short, we wrote the present note with the
description of the algorithm and the proof of its correctness.

Let each $\B x_i$ be represented by its Cartesian coordinates, viewed
as a column vector: $\B x_i=(x_{i1},x_{i2},x_{i3})^T$. Clearly, the
existence a supporting plane is equivalent to the existence of a
linear function $f(\B y)=\sum_{j=1}^3 f_iy_i$ with $f(\B x_i)>0$ for
all $1\le i\le 4$.

It is obvious how to extend the above notions to the general case of
$k$ vectors in $\I R^d$. Here is a criterion that we will use. 

\begin{lemma}\label{lm:\pl:weights} Let $\B x_1,\dots,\B x_k\in\I R^d$. A supporting
$(d-1)$-dimensional hyperplane does not exist if and only if there
are weights $w_i\ge 0$, $1\le i\le k$, not all zero, with
$\sum_{i=1}^k w_i\B x_i=\B 0$.
 \end{lemma}
 \smallskip{\it Proof.}  If a supporting plane does not exist, then the system of linear
inequalities $\sum_{j=1}^d f_ix_{ij}\ge 1$ has no solution in
the $f$'s. Now the existence of weights is precisely the conclusion of
the Farkas Lemma. (For example, apply Proposition~1.7
from~\cite{ziegler:lp} with $A_{ij}=-x_{ij}$ and $z_i=-1$, $1\le i\le
k$, $1\le j\le d$.)

On the other hand, if we have weights, then for any $f$ we have
$\sum_{i=1}^k w_i f(\B x_i) =f(\sum_{i=1}^k w_i\B x_i)=\B 0$, so $f(\B
x_i)\le 0$ for some $i$.\qed \medskip

Return to our $\B x_1,\dots,\B x_4\in\I R^3$. First, we handle the
case when some three vectors are coplanar.

\begin{theorem}\label{th:\pl:coplanar} Let $\B x_1$, $\B x_2$ and $\B x_3$ lie on a plane
$P$. There exists a supporting plane if and only if the set $\{\B
x_1,\dots,\B x_4\}\cap P$ admits a supporting line in $P$.\end{theorem}
 \smallskip{\it Proof.}  If $\B x_4\in P$, then both weight conditions of
Lemma~\ref{lm:\pl:weights} are identical. If $\B x_4\not\in P$, then a
supporting plane does not exist iff $\sum_{i=1}^4 w_i\B x_i=\B 0$ for
some weights, which is the case iff $\sum_{i=1}^3 w_i\B x_i=\B 0$ for
some weights because $\B x_4$ does not lie in the linear span of $\B
x_1,\B x_2,\B x_3$.\qed \medskip

Testing in a plane is probably done easiest through the following
observation: a supporting line for vectors $\B x_1,\dots,\B x_k\in \I
R^2$ exists iff there is a relabelling of subscripts such that
$\angle(\B x_1,\B x_k)=\sum_{i=1}^{k-1}\angle(\B x_i,\B x_{i+1})$ and
$\angle(\B x_1,\B x_k)<\pi$, where $0\le \angle(\B x,\B y)\le \pi$ is
the angle between vectors $\B x$ and $\B y$.

It remains to consider the case when no three vectors are
coplanar. Let $P(\B x_i,\B x_j)$ be the plane passing through $\B x_i$
and $\B x_j$. Here is a corollary to Lemma~\ref{lm:\pl:weights}.

\begin{corollary}\label{cr:\pl:Pij} A supporting plane exists if and only if, for some $1\le i<
j\le 4$, the remaining two vectors lie in the same open halfspace with
respect to the plane $P(\B x_i,\B x_j)$.\end{corollary}
 \smallskip{\it Proof.}  Suppose we have a supporting plane $P$. We can rotate $P$ until
it hits some $\B x_i$. If $P$ does not contain other $\B x$'s, rotate
it further, along the axis $\B 0 \B x_i$, until it hits for the first
time some $\B x_j$. Now, $P(\B x_i,\B x_j)$ is clearly the required
plane. 

Conversely, suppose that a supporting plane does not exist, that is,
by Lemma~\ref{lm:\pl:weights} we can find weights. Let $f=0$ be the equation
of the plane $P(\B x_1,\B x_2)$. We have $f(\B x_3)\not=0$ and $f(\B
x_4)\not=0$ and the equality $\sum_{h=1}^4 w_h f(\B x_h)=0$ implies
that $f(\B x_3)\cdot f(\B x_4)<0$, as required.\qed \medskip

So, the algorithm would be to take all six pairs $1\le i<j \le 4$,
write an equation of the plane $P(\B x_i,\B x_j)$, which is
 $$
 f(\B y)=\det(\B x_i,\B x_j,\B y)=0,
 $$
 and check the signs of $f$ on the remaining two vectors. However,
this procedure is not very economical as the determinant of a matrix,
say $(\B x_1,\B x_2,\B x_3)$, is computed essentially three
times. Here is a better algorithm.

\begin{theorem}\label{th:\pl:algorithm} Compute the signs of
 $$ \det(\B x_1,\B x_2,\B x_3),\ \det(\B x_1,\B x_4,\B x_2),\ \det(\B
x_3,\B x_2,\B x_4),\ \det(\B x_1,\B x_3,\B x_4).$$
 A supporting plane does not exist if and only if all four signs are
the same.\end{theorem}
 \smallskip{\it Proof.}  Suppose that a supporting plane does not exist. Assume
$\det(\B x_1,\B x_2,\B x_3)>0$. Then applying Corollary~\rcr{Pij} to
$P(\B x_1,\B x_2)$ we obtain that $\det(\B x_1,\B x_2,\B x_4)<0$, that
is, $\det(\B x_1,\B x_4,\B x_2)>0$ as required. Similarly, $\det(\B
x_3,\B x_2,\B x_4)>0$ and $\det(\B x_1,\B x_3,\B x_4)>0$.

The converse also follows by an easy application of
Corollary~\rcr{Pij}: e.g.\ $\B x_3$ and $\B x_4$ are separated by
$P(\B x_1,\B x_2)$ because
 $$\det (\B x_1,\B x_2,\B x_3)\cdot \det(\B x_1,\B x_2,\B x_4) = -
\det (\B x_1,\B x_2,\B x_3)\cdot \det(\B x_1,\B x_4,\B x_2) < 0.\qed
 $$
 
\providecommand{\bysame}{\leavevmode\hbox to3em{\hrulefill}\thinspace}

\end{document}